\input amstex
\documentstyle {amsppt}
\magnification=1200
\nologo
\NoBlackBoxes
\document
\pageheight{18cm}
\centerline{\bf STABLE MAPS OF GENUS ZERO}

\smallskip

\centerline{\bf TO FLAG SPACES}

\medskip

\centerline{Yu.I.Manin}

\medskip

\centerline{Max--Planck--Institut f\"ur Mathematik, Bonn}

\bigskip

\hfill{\it Dedicated to J.~Moser}

\medskip

\centerline{\bf \S 0. Introduction}

\medskip

Topological quantum field theory led recently to a spectacular
progress in numerical algebraic geometry. It was shown 
that generating functions of certain charactertistic numbers
of modular spaces of stable algebraic curves with labelled points
satisfy remarkable differential equations of KP--type
(E. Witten, M. Kontsevich). In a later series of developments,
this was generalized, partly conjecturally, to the spaces of maps of curves
into algebraic varieties leading to the Mirror Conjecture
and the construction of quantum cohomology.

\smallskip

The key technical notion in the context of algebraic
geometry is that of a {\it stable map} introduced
by M. Kontsevich (cf. [K] and [BM]) following the earlier
work by M. Gromov in symplectic geometry. It provides a natural
compactification of spaces of maps, in the same way as
stable curves compactify moduli spaces.

\smallskip

We will be working over a ground field. Let $W$ be an
algebraic variety.

\smallskip

\proclaim{\quad 0.1. Definition} A stable map (to $W$) is a
structure $(C;x_1,\dots ,x_n;f)$ consisting of the following data.

\smallskip

a). $(C;x_1,\dots ,x_n)$ is a connected complete reduced curve
with $n\ge 0$ labelled pairwise distinct non--singular points
$x_i$ and at most ordinary double singular points.

\smallskip

b). $f:\ C\to W$ is a morphism having no non--trivial first
order infinitesimal automorphisms identical on $W$ and
$x_i$'s (stability). This means that every irreducible component
of $C$ of genus zero (resp. 1) has at least three (resp. one)
special points (inverse images of singular and labelled points)
on its normalization.

\endproclaim

\smallskip

A family of stable maps parametrized by a noetherian scheme $S$
is a structure consisting of a flat proper morphism
$\pi :\ \Cal{C}\to S$, $n$ sections $x_i:\ S\to \Cal{C}$,
and a morphism $f:\ \Cal{C}\to W$ whose restriction to each geometric
fiber of $\pi$ is a stable map in the sense of the previous
Definition.

\smallskip

Families of stable maps form an algebraic stack in the sense
of [DM].

\smallskip

For  fixed $n\ge 0,\ g\ge 0$ and an algebraic homology class
of dimension two $\beta$, denote by $\overline{\Cal{M}}_{g,n}(W,\beta )$
the substack of maps for which $g$ is the arithmetic genus of $C$
and $\beta =f_*([C]).$ For the proof of the following theorem see
[K], [BM]  and [FP].

\smallskip

\proclaim{\quad 0.2. Theorem} a). If $W$ is projective, then
$\overline{\Cal{M}}_{g,n}(W,\beta )$ is a proper separated algebraic stack
of finite type.

\smallskip

b). Assume that $W$ is convex in the following sense:
$H^1(C,f^*(\Cal{T}_W))=0$ for any stable map $f:\ C\to W$
of genus zero. Then the stacks $\overline{\Cal{M}}_{0,n}(W,\beta )$
are smooth.

\smallskip

Denote by $\Cal{M}_{0,n}(W,\beta )$ the big cell of this stack
over which $C$ is smooth. The complement of this cell is
a divisor with normal crossings.
\endproclaim

In addition, the spaces of geometric points of these stacks
are represented by algebraic schemes, crude moduli spaces
of stable maps, $\overline{M}_{g,n}(W,\beta )$ and
$M_{g,n}(W,\beta )$ respectively. If $(C;x_1,\dots ,x_n;f)$
is a stable map,
we denote by $[(C;x_1,\dots ,x_n;f)]$, or simply $[f]$,
the corresponding point. Two maps $[f^{(i)}],\ i=1,2,$ define the same
point iff there is an isomorphism $g:\ C^{(1)}\to C^{(2)}$
such that $g(x_i^{(1)})=x_i^{(2)}$ for all $i$ and
$g\circ f^{(2)}=f^{(1)}.$

\smallskip

This paper is a sequel to [M].

\smallskip

We will consider mostly the case $g=0$, $W=G/P$
(generalized flag spaces). Although flag spaces are
convex, the respective spaces of stable maps are not
smooth but only orbifolds in general. Our goal is to calculate
their virtual Poincar\'e polynomials, or rather an
appropriate generating function for these polynomials.
This calculation generalizes the one made for $\overline{M}_{0,n}$
in [M] and [G]. Remarkably,
up to a change of variables, this function satisfies the same
universal differential equation as that in [M], the dependence
on $W$ being reflected only in the initial condition which
involves the Eisenstein series of $W=G/P$: see 2.3 below.

\smallskip

{\it Acknowledgement.} I am grateful to Nikita Nekrasov who alerted me
to several misprints and errors in the first draft of this paper.

\bigskip

\centerline{\bf \S 1. Stratification of the space of stable maps}

\bigskip

{\bf 1.1. Virtual Euler--Poincar\'e maps.} Let $\Cal{R}$
be a commutative associative $\bold{Q}$--algebra, $Var$
the category of algebraic varieties, not necessarily complete,
smooth, or irreducible, over a fixed ground field.

\smallskip

We will call a map $Ob\ Var\to \Cal{R}:\ V\mapsto [V]$
a {\it virtual Euler--Poincar\'e map} if the following conditions
are satisfied:

\smallskip

i). $[V]$ depends only on the isomorphism class of $V$.

\smallskip

ii). Additivity: if $V=\coprod_iV_i$ is a finite disjoint
union of locally closed subsets (``strata''),
then $[V]=\sum_i[V_i].$

\smallskip

iii). Multiplicativity: $[\prod_jV_j]=\prod_j[V_j].$

\smallskip

From ii) and iii) it follows that if $V\to B$ is a Zariski
locally trivial fiber space with fiber $F,$ then
$[V]=[B][F].$

\smallskip

If we work over $\bold{C},$ and $\Cal{R}=\bold{Q}[q]$ (resp. $\Cal{R}=
\bold{Q}[r,s]$) the main examples are virtual
Poincar\'e (resp. Hodge) polynomials:
$$
P_V(q):=\sum_{i,j}(-1)^{i+j}\roman{dim}(\roman{gr}^j_{W}
H^i_c(V))q^j, \eqno{(1.1)}
$$
$$
H_V(r,s):=\sum_{i,j,k}(-1)^{i+j+k}h^{j,k}
(H^i_c(V))r^js^k, \eqno{(1.2)}
$$
and also the virtual Euler characteristics $\chi (V)=P_V(-1).$
These properties were for the first time systematically
used by Danilov and Khovanskii in the toric geometry.

\smallskip

Recently C. Soul\'e and H. Gillet ([GS]) established that the map
$V\mapsto$ {\it class of the motive} $h^*(V)$ {\it of}\ $V$ {\it in the}
$K_0$--{\it ring of Grothendieck's motives} extends from projective
smooth varieties to $Var$ and becomes an Euler--Poincar\'e
map. For a somewhat different construction see [GuNA]. 

\smallskip

On the subcategory generated by the (dual) Tate motive one can
identify the Gillet--Soul\'e map with virtual
Poincar\'e map putting $h^*(\bold{P}^1)=q^2+1$,
and calculate the latter via point count over $\bold{F}_{q^2}.$
Most of our calculations are in fact restricted to this
situation.

\smallskip

We will have also to localize $\Cal{R}$, most notably by
inverting $[PGL(2)]=q^2(q^4-1).$ Since the Euler characteristic
of this manifold vanishes, we have to apply a limiting procedure
$q^2\to 1$ producing logarithms and infinite ramification
in our formulas for generating functions of $\chi$,
as already happened in [M].

\medskip

{\bf 1.2. Poincar\'e polynomials of $\underline{Map}_{\beta}
(\bold{P}^1,W).$} If $W$ is a flag manifold, the scheme
of maps of $\bold{P}^1$ to $W$ landing in $\beta$
is a smooth (generally non--complete) manifold. Its
Gillet--Soul\'e motive lies in the Tate's subring,
and the generating function
$$
E(W,z):=\sum_{\beta\in B}[\underline{Map}_{\beta}(\bold{P}^1,W)]z^{\beta}
\eqno{(1.3)}
$$
is rational in $z$. Here $B\subset N$ means the effective subsemigroup 
of the lattice $N$ of algebraic homology classes of curves, 
and $z^{\beta}$ are formal monomials, elements of the dual 
lattice $M$ written multiplicatively.

\smallskip

In fact, in the $\bold{F}_{q^2}$--avatar, (1.3) essentially
coincides with an appropriate Eisenstein series discussed
e.g. in [H]. If we introduce a complex
vector variable $s\in M_{\bold{C}}$ and replace $z^{\beta}$ by 
$q^{-2(\beta .s)}$, the poles of Eisenstein series will
lie on the complexified walls of the ample cone $B^t_{\bold{R}}$ in $N$
shifted by $-K_W,$
and reflections with respect to the walls will generate the
functional equations.

\medskip

{\bf 1.2.1. Example.} We have
$$
E(\bold{P}^n,z)=\frac{1-q^{2n+2}}{1-q^2}\frac{1-q^2z}{1-q^{2n+2}z}.
\eqno{(1.4)}
$$
In fact, denote by $N_d$ the number of non--zero $(n+1)$--tuples
$(f_0(t_0,t_1),\dots ,f_n(t_0,t_1))$ of coprime forms of degree
$d$ over $\bold{F}_{q^2}$ divided by $|\bold{F}_{q^2}^*|=q^2-1.$
Since the number of all such $(n+1)$--tuples is
$q^{2(n+1)(d+1)}-1$ and the  map $((f_0,\dots ,f_n),g)\mapsto
(f_0g,\dots ,f_ng)$ with fixed $g\ne 0$ and degrees of $f_i,g$
has all fibers of cardinality $q^2-1$, we have
$$
\sum_{k=0}^{d}N_{d-k}(q^{2(k+1)}-1)=q^{2(k+1)(d+1)}-1,
$$
from which (1.4) easily follows.

\smallskip

Identifying $M=\roman{Pic}\ {\bold{P}^n}$ with $\bold{Z}$
via $\Cal{O}(1)\mapsto 1$ and putting $z=q^{-2s}$ one sees
that the real pole of (1.4) is $s=n+1=-K_{\bold{P}^n}.$

\medskip

{\bf 1.3. Big cells of stable map spaces.} We will say that
a stable map $(C;x_1,\dots ,x_n;f)$ to $W$ belongs to the
big cell $M_{0,k}(W,\beta )$ if $C\cong \bold{P}^1$, and 
$f_*([C])=\beta $. We will explain how these big cells 
are related to $\underline{Map}_{\beta}(\bold{P}^1,W).$

\smallskip

a). $\beta =0.$ In this case, $k\ge 3$, and the map
$$
(\bold{P}^1;x_1,\dots ,x_k;f)\mapsto (f(x_1),
[(\bold{P}^1;x_1,\dots ,x_k)])
$$
induces the
identification
$$
M_{0,k}(W,0)\cong W\times M_{0,k}.   \eqno{(1.5)}
$$
In particular, when $W$ is a point, we get simply $M_{0,k}.$

\smallskip

b). $\beta \ne 0.$ In this case, a choice of three points $x_1,x_2,x_3$
in $\bold{P}^1$ defines an isomorphism
$$
M_{0,3}(W,\beta)\cong \underline{Map}_{\beta}(\bold{P}^1,W).
\eqno{(1.6)}
$$
Let $G^{(i)}\subset PGL(2)$ be the subgroup fixing $x_1,\dots ,x_i,\ 
i=0,\dots ,3$. Then $G^{(i)}$ freely acts upon $M_{0,3}(W,\beta )$,
and we have
$$
M_{0,i}(W,\beta )\cong \underline{Map}_{\beta}(\bold{P}^1,W)/G^{(i)},\
i\le 3.   \eqno{(1.7)}
$$
Finally, for $k\ge 4$ the forgetful map $M_{0,k}(W,\beta )\to
M_{0,3}(W,\beta )$ identifies $M_{0,k}(W,\beta )$ with a locally
trivial fibration over $\underline{Map}_{\beta}(\bold{P}^1,W)$
with fiber 
$$
M_{0,k}=((\bold{P}^1)^k\setminus\cap_{i<j}\Delta_{ij})/
PGL(2).
$$

\smallskip

To summarize, we have the following formula for the virtual 
Euler--Poincar\'e class of $M_{0,k}(W,\beta)$ valid for
$\beta =0$ as well:

\smallskip

\proclaim{\quad 1.3.1. Proposition} We have
$$
[M_{0,k}(W,\beta )]=[\underline{Map}_{\beta}(\bold{P}^1,W)]
{{[\bold{P}^1]}\choose{k}}k!\frac{1}{[PGL(2)]}.
\eqno{(1.8)}
$$
\endproclaim

\medskip

In fact,
$$
[M_{0,k}]={{[\bold{P}^1]-3}\choose{k-3}}(k-3)!,\qquad
[PGL(2)]={{[\bold{P}^1]}\choose{3}}3!.
$$
We will now construct strata of the space of stable maps numbered
by marked trees. Their virtual Euler--Poincar\'e classes
will be expressed via products of those of big cells. 

\medskip

{\bf 1.4. Trees.} In this paper, a {\it tree} $\tau$ is a
finite connected simply connected CW--complex of dimension
$1$ or $0$ (one--vertex tree). A {\it flag} of a tree is
a pair consisting of a vertex and an adjoining edge.
The {\it valency} $|v|$ of a vertex $v$ is the number
of flags containing $v$. The sets of vertices (resp. edges, flags)
of $\tau$ are denoted $V_{\tau}$ (resp. $E_{\tau},F_{\tau}$).

\smallskip

For a set $S$, $\Cal{P}(S)$ denotes the set of its subsets.
As above, $B$ is the semigroup of classes of effective
algebraic curves in a flag manifold $W$.

\smallskip

\proclaim{\quad 1.4.1. Definition} A $(k,W)$--marking of
a tree $\tau$ is a map
$$
\mu :\ V_{\tau}\to B\times \Cal{P}({1,\dots ,k}):\ 
v\mapsto (\beta_v,S_v),
$$
satisfying the following conditions:

\smallskip

a). The family $\{ S_v|v\in V_{\tau}, S_v\ne \emptyset\}$
forms a partition of $\{ 1,\dots ,k\}$ into pairwise
disjoint subsets. (For $k=0$, we interpret $\Cal{P}$
as one--element set, and this condition as empty).

\smallskip

b). If $\beta_v=0,$ then $|v|+|S_v|\ge 3.$
\endproclaim

\smallskip

There is an obvious notion of isomorphism of marked trees.

\medskip

{\bf 1.4.1. Type of a stable map.} Let $(C;x_1,\dots ,x_k;f)$
be a stable map of genus $0$ to $W.$ Its {\it combinatorial
type} is, by definition, the isomorphism class of the dual tree
of $C$ with the marking depending on $f.$ We can describe it
as follows. First of all, put
$$
V_{\tau}=\{ irreducible\ components\ C_v\ of\ C\},
$$
$$
E_{\tau}=\{ intersection\ points\ of\ irreducible\ components\}.
$$
Thus a flag of $\tau$ is a pair {\it (intersection point
of two components, one of these components)}. Furthermore,
$$
\beta_v:=f_*([C_v]),\ S_v=\{ i\,|\,x_i\in C_v\}.
$$
Obviously, $\beta_v=0$ means that $f(C_v)$ is a point so that
condition 1.4.1 b) means stability.

\medskip

{\bf 1.5. Maps of fixed type.} Consider a combinatorial type
$(\tau ,\mu ).$ For any vertex $v\in V_{\tau}$, construct
the stack $M_{0,S_v}(W,\beta_v)$ parametrizing stable
maps of $\bold{P}^1$ with points labelled by $S_v$,
of class $\beta_v$. For any $t\in S_v$, there is a canonical
evaluation map $ev_t:\,M_{0,S_v}(W,\beta_v) \to W:\
[f]\mapsto f(x_t)$ where $x_t$ is the point marked by $t$.
Put
$$
M_{0,k}(W,\beta ,\tau ,\mu ) = \prod_W M_{0,S_v}(W,\beta_v),
$$
where $\prod_W$ means the partial fibered product over $W$
which leaves in the total product only those
families of stable maps $(f_v\,|\,v\in V_{\tau})$
for which $f_v(s)=f_w(t)$ whenever $s,t$ are halves of the same
edge of $\tau$.

\smallskip

\proclaim{\quad 1.5.1. Theorem} Let $W$ be a generalized
flag variety. Put
$$
N(W,\beta ):=\frac{[M_{0,3}(W,\beta)]}{[W][PGL(2)]}=
\frac{[\underline{Map}_{\beta}(\bold{P}^1,W)]}{[W][PGL(2)]}.
\eqno{(1.9)}
$$
Then
$$
[M_{0,k}(W,\beta ,\tau ,\mu )]=
$$
$$
[W]\prod_{v\in V_{\tau}}\varepsilon (\beta_v,|v|+k_v)
N(W,\beta_v){{[\bold{P}^1]}\choose{|v|+k_v}}(|v|+k_v)!
\eqno{(1.10)}
$$
where $\mu :\ v\mapsto (\beta_v,S_v),k_v=|S_v|,\beta =\sum\beta_v,$
and $\varepsilon (\beta ,n)=0$ for $\beta =0,n=0,1,2,$ and $1$
otherwise.
\endproclaim

\smallskip

{\bf Proof (sketch).} For $\beta =0,$ we have a projection
$M_{0,k}(W,0,\tau ,\mu )\to W:\ (C;x_1,\dots ,x_k;f)\mapsto f(C).$
Its fiber is the stratum of the moduli space of stable curves
with labelled points of the given combinatorial type $(\tau ,\mu )$;
notice that $\mu$ now is just a map $v\mapsto S_v$. This fiber
is isomorphic to $\prod_{v\in V_{\tau}}M_{0,|v|}.$ The fibration
is locally trivial over $[W].$ Hence
$$
[M_{0,k}(W,0,\tau ,\mu )]=[W]\prod_{v\in V_{\tau}}
{{[\bold{P}^1]-3}\choose{|v+k_v|-3}}(|v+k_v|-3)!.
$$
One easily sees that this coincides with (1.10).

\smallskip

For $\beta\ne 0,\tau$ one--vertex that is, $C=\bold{P}^1,$
(1.10) follows from (1.8).

\smallskip

Finally, in the remaining cases we represent a stable map $f:C\to W$ as
a vector of stable maps of irreducible components $f_v:C_v\to W$.
If we fix one component, the rest will be constrained by
incidence conditions. This accounts for the necessity to divide
by $[W]$ in each $v$--factor; we incorporated this division
in the definition (1.9) of $N(W,\beta ).$ The exterior
factor $[W]$ in (1.10) appears because one component can be moved
freely.

\bigskip

\centerline{\bf \S 2. Generating function}
\centerline{\bf of the space of stable maps}

\bigskip

{\bf 2.1. The problem.} In this section, we set to calculate
the following series in $\Cal{R}[[t,z]]$ (we remind
that monomials in $z$ belong to a semigroup ring):
$$
\Phi_W(t,z):=\sum_{\tau/(iso)}\frac{1}{|\roman{Aut}\,\tau |}
\sum_{{V_{\tau}\to B\times N}\atop{v\mapsto (\beta_v,k_v)}}[W]\times
$$
$$
\prod_{v\in V_{\tau}}\frac{t^{k_v}}{k_v!}z^{\beta_v}
\varepsilon(\beta_v,|v|+k_v)N(W,\beta_v)
{{[\bold{P}^1]}\choose{|v|+k_v}}(|v|+k_v)!
\eqno{(2.1)}
$$
We retain the notation of \S 1; in particular, $\varepsilon$
is introduced in order to exclude inadmissible markings.

\smallskip

To motivate this definition, consider the following more natural
generating function:
$$
\widetilde{\Phi}_W(t,z):=\sum_{\beta ,k}
[\overline{M}_{0,k}(W,\beta )]\frac{t^k}{k!}z^{\beta }.
\eqno{(2.2)}
$$
Assume in addition that $[\overline{M}_{0,k}(W,\beta )]$
in (2.2) means an ``orbifold'' virtual class of 
$[\overline{M}_{0,k}(W,\beta )]$ which is defined on the
category of small algebraic stacks and besides the
usual additivity and multiplicativity properties
postulated in 1.1 satisfies the following condition:
{\it if} $V$ {\it is smooth,} $G$ {\it a finite group freely
acting upon} $V$, {\it then} $[V/G]=\dfrac{[V]}{|G|}.$

\smallskip

In this case (2.2) coincides with (2.1). To see this, one has
to consider the stratified covering of $\overline{M}_{0,k}(W,\beta )$
by $M_{0,k}(W,\beta ,\tau ,\mu )$ where $(\tau ,\mu )$ runs
over all admissible marked trees with $\sum_v\beta_v=\beta ,
\sum_vk_v=k,$ and $M_{0,k}(W,\beta .\tau ,\mu )$ are moduli spaces
of rigidified stable maps. The claim follows from the fact
that $Aut\ \tau$ acts freely on $\coprod_{\mu}M_{0,k}(W,
\beta ,\tau ,\mu ).$

\smallskip

\proclaim{\quad 2.2. Theorem} Denote by $\varphi^0=\varphi_W(t,z)$
the unique root in $\Cal{R}[[t,z]]$ of the following equation:
$$
\frac{[\bold{P}^1]}{[PGL(2)][W]}E(W,z)(1+t+\varphi^0)^{[\bold{P}^1]-1}=
$$
$$
\varphi^0\frac{[\bold{P}^1]-1}{[\bold{P}^1]-2}
+\frac{t}{[\bold{P}^1]-2}
+\frac{1}{([\bold{P}^1]-1)([\bold{P}^1]-2)}
\eqno{(2.3)}
$$
satisfying the condition
$$
\varphi^0\equiv\sum_{\beta\in ind}N(W,\beta )z^{\beta}
\roman{mod}\ (z,t)^2,
$$
where the summation is taken over indecomposable elements
of $B$. Then we have
$$
\Phi_W(t,z)=-\frac{[\bold{P}^1]-1}{2[\bold{P}^1]}(\varphi^0)^2
+\frac{\varphi^0}{[\bold{P}^1]}
-\frac{t^2}{2[\bold{P}^1]}
\eqno{(2.4)}
$$
$$
{\frac{\partial \Phi_W}{\partial t}}(t,z)=\varphi^0.
\eqno{(2.5)}
$$
\endproclaim
\smallskip

{\bf Proof.} We first rearrange the inner sum and product in (2.1):
$$
\Phi_W(t,z)=\sum_{\tau /(iso)}\frac{[W]}{|Aut\ \tau|}
\prod_{v\in V_{\tau}}\sum_{\beta ,k\ge 0}
\frac{t^k}{k!}z^{\beta}\times
$$
$$
\varepsilon(\beta ,|v|+k)N(W,\beta){{[\bold{P}^1]}\choose{|v|+k}}
(|v|+k)!
\eqno{(2.6)}
$$
Furthermore, for a fixed $v\in V_{\tau}$,
$$
\sum_{\beta ,k}\frac{t^k}{k!}z^{\beta}\varepsilon (\beta ,|v|+k)
N(W,\beta ){{[\bold{P}^1]}\choose{|v|+k}}(|v|+k)!=
$$
$$
\left( \sum_{\beta}N(W,\beta )z^{\beta}\right)
\left(\sum_{k\ge 0}\frac{t^k}{k!}{{[\bold{P}^1]}\choose{|v|+k}}
(|v|+k)!\right) -\varepsilon_{|v|}=
$$
$$
\frac{E(W,z)}{[PGL(2)][W]}\sum_{k\ge 0}\frac{t^k}{k!}
{{[\bold{P}^1]}\choose{|v|+k}}(|v|+k)!-\varepsilon_{|v|},
\eqno{(2.7)}
$$
where $\varepsilon_n=0$ for $n\ge 3,$ and
$$
\varepsilon_n=\sum_{k=0}^{2-n}N(W,0)\frac{t^k}{k!}
{{[\bold{P}^1]}\choose{n+k}}(n+k)!\qquad \roman{for}\ n\le 2.
\eqno{(2.8)}
$$
To calculate the resulting sum over marked trees we can now
apply the formalism of [K], [M]. Introduce one more formal
variable $\varphi$ and consider the formal potential
$$
S(\varphi ):=-\frac{\varphi^2}{2}+\sum_{n=0}^{\infty}
\frac{C_n}{n!}\varphi^n,
\eqno{(2.9)}
$$
$$
C_n:=\frac{E(W,z)}{[PGL(2)][W]}\sum_{k=0}^{\infty}
\frac{t^k}{k!}(n+k)!{{[\bold{P}^1]}\choose{n+k}}-\varepsilon_n.
\eqno{(2.10)}
$$
We have
$$
S(\varphi )=-\frac{\varphi^2}{2}+\frac{E(W,z)}{[PGL(2)][W]}
\sum_{n,k=0}^{\infty}\varphi^nt^k
{{n+k}\choose{k}}{{[\bold{P}^1]}\choose{n+k}}-
\sum_{n=0}^{2}\frac{\varepsilon_n}{n!}\varphi^n=
$$
$$
\frac{E(W,z)}{[PGL(2)][W]}(1+t+\varphi )^{[\bold{P}^1]}-
\varphi^2\frac{[\bold{P}^1]-1}{2([\bold{P}^1]-2)}
$$
$$
-\varphi \left(\frac{1}{([\bold{P}^1]-1)([\bold{P}^1]-2)}
+\frac{t}{[\bold{P}^1]-2}\right)
$$
$$
-\left(\frac{1}{[\bold{P}^1]([\bold{P}^1]-1)([\bold{P}^1]-2)}
+\frac{t}{([\bold{P}^1]-1)([\bold{P}^1]-2)}
+\frac{t^2}{2([\bold{P}^1]-2)}\right) .
\eqno{(2.11)}
$$
According to a general formula of perturbation theory (cf. [K], [M])
we have
$$
\frac{1}{[W]}\,\Phi_W(t,z)=S^{crit}:=S(\varphi^0),
$$
where $\varphi^0$ is an appropriate formal solution of
$\dfrac{\partial}{\partial \varphi}S(\varphi )=0$ which we will
now identify. Differentiating (2.11) in $\varphi$ we obtain:
$$
\frac{E(W,z)}{[PGL(2)][W]}[\bold{P}^1](1+t+\varphi^0)^{[\bold{P}^1]-1}=
$$
$$
\varphi^0\frac{[\bold{P}^1]-1}{[\bold{P}^1]-2}
+\frac{t}{[\bold{P}^1]-2}
+\frac{1}{([\bold{P}^1]-1)([\bold{P}^1]-2)}
$$
which is (2.3).

\smallskip

Rewriting this equation as
$$
\varphi^0=C_1+C_2\varphi^0+C_3\frac{(\varphi^0)^2}{2}+\dots
$$
and taking into account that $C_1$ modulo $(t,z)^2$ starts
with terms linear in $z$ we see that there is a unique
solution $\varphi^0\in\Cal{R}[[t,z]]$ such that
$\varphi^0\equiv C_1\roman{mod}\ (t,z)^2.$ In view
of (2.10), this coincides with the congruence in the
statement of the Theorem.

\smallskip

It remains to calculate $S(\varphi^0).$ Multiplying (2.3)
by $\dfrac{1+t+\varphi^0}{[\bold{P}^1]}$
and simplifying, we get (2.4).

\smallskip

Finally, derivating (2.3) and (2.4) in $t$, we can obtain
(2.5): we reproduce this calculation below in simplified
notation.

\medskip

{\bf 2.3. Comments and supplements.} We will now consider
the case $\Cal{R}=\bold{Q}[q], [V]=P_V(q)$
(the virtual Poincar\'e polynomial of $V$). Equations
(2.3) and (2.4) take form
$$
\frac{1}{q^2(q^2-1)P_W(q)}E(W,z)(1+t+\varphi^0)^{q^2}=
\varphi^0\frac{q^2}{q^2-1}+\frac{t}{q^2-1}+\frac{1}{q^2(q^2-1)},
\eqno{(2.12)}
$$
$$
\frac{1}{[P_W(q)]}\,\Phi_W(t,z)=-\frac{q^2}{2(q^2+1)}(\varphi^0)^2
+\frac{1}{q^2+1}\varphi^0-\frac{1}{2(q^2+1)}t^2.
\eqno{(2.13)}
$$

\smallskip

a). {\it A differential equation for} $\varphi^0.$
Derivating (2.12) in $t$ and multiplying the result
by $\dfrac{1+t+\varphi^0}{q^2}$ we get the differential
equation
$$
(1-q^2\varphi^0)\varphi^0_t=(q^2+1)\varphi^0+t.
\eqno{(2.14)}
$$
Up to a simple variable change $\varphi^0=\psi -t$,
this is the same equation as (0.7) (and (0.15)) in [M]:
$$
(1+q^2t-q^2\psi )\psi_t=1+\psi .
\eqno{(2.15)}
$$
Its universality is remarkable: in (2.14) there is no dependence
on $W$ and $z$ (below it will be encoded in the choice of 
constant in a general solution of (2.14)), and in [M]
it emerged also in a calculation of the generating function for
arbitrary configuration spaces $X[n]$.

\smallskip

b). {\it Euler characteristics:} $q^2=1$.
A formal substitution $q^2=1$ into (2.11) is impossible,
however a limiting procedure gives:
$$
(1+t+\varphi^0)\roman{log}\,(1+t+\varphi^0)=
2\varphi^0 +t- X(1+t+\varphi^0),
$$
$$
 X =\left. 
\frac{1}{q^2-1}\,\left(\frac{E(W,z)}{P_W(q)}-1\right)\right|_{q^2=1}.
\eqno(2.16)
$$
and $\Phi_W(t,z)$ becomes
$$
\frac{1}{[P_W(q)]}\,\Phi_W(t,z)=-\frac{1}{4}(\varphi^0)^2+\frac{1}{2}\varphi^0
-\frac{1}{4}t^2 . \eqno({2.17)}
$$

\smallskip

c). {\it General solution of (2.14).} We can apply to (2.14)
Proposition 1.5.1 in [M]. Put
$$
x=t+\frac{q^2+1}{q^2},\ y=q^2\varphi^0-1,\ w=\frac{y}{x}.
$$
The general solution of (2.14) in implicit form is
given by
$$
Cx=(w+1)^{\frac{1}{q^2-1}}(w+q^2)^{\frac{q^2}{1-q^2}}.
\eqno{(2.18)}
$$
This makes evident the ramification structure of $\varphi$
for $q^2\ne 1.$ The constant $C$ in our case is a function of $z$
which can be calculated by considering (2.18) at $t=0$ that is,
$x=\dfrac{q^2+1}{q^2}.$ The function $w(z)=\dfrac{q^4\varphi^0(0,z)-q^2}{q^2+1}$
is found from (2.12) at $t=0.$

\smallskip

d). As in [M], the derivative of $\Phi_W$ in $t$ is simpler
than $\Phi_W$ itself. Indeed, from (2.13) and (2.14) we have:
$$
\frac{\partial \Phi_W}{\partial t}=
-\frac{q^2}{q^2+1}\varphi^0\varphi^0_t
+\frac{1}{q^2+1}\varphi^0_t-\frac{t}{q^2+1}
$$
$$
=-\frac{1}{q^2+1}\left(\varphi^0_t-(q^2+1)\varphi^0-t\right)
+\frac{1}{q^2+1}\varphi^0_t-\frac{t}{q^2+1}=\varphi^0.
$$

\newpage

\centerline{\bf References}

\medskip

[BM] K.~Behrend, Yu.~Manin. {\it Stacks of stable maps and
Gromov--Witten invariants.} Duke Math. J., 85:1 (1996), 1--60.

\medskip

[DM] P.~Deligne, D.~Mumford. {\it The irreducibility
of the space of curves of a given genus.} Publ. Math. IHES,
36 (1969), 75--109.

\medskip

[FP] W.~Fulton, R.~Pandharipande. {\it Notes on stable
maps and quantum cohomology.} Preprint alg--geom/9608011.

\medskip

[G] E.~Getzler. {\it Operads and moduli spaces of genus zero
Riemann surfaces.} In: The Moduli Space of Curves, ed. by
R\. Dijkgraaf, C\. Faber, G\. van der Geer, Progress in Math\.
vol\. 129, Birkh\"auser, 1995, 199--230.

\medskip

[GiS] H.~Gillet, C.~Soul\'e. {\it Descent, motives and K--theory.}
J. reine u. angew. Math., 478 (1996), 127--176.

\medskip

[GuNA] F.~Guill\'en, V.~Navarro Aznar. {\it Un crit\`ere d'extension
d'un foncteur d\'efini sur sch\'emas lisses.} Preprint, 1995.

\medskip

[H] G.~Harder. {\it Chevalley groups over function fields and
automorphic forms.}
Annals of Math.,  100(1974), 249--306.

\medskip

[K] M.~Kontsevich. {\it Enumeration of rational curves via torus actions.}
In: The Moduli Space of Curves, ed. by
R\. Dijkgraaf, C\. Faber, G\. van der Geer, Progress in Math\.
vol\. 129, Birkh\"auser, 1995, 335--368.

\medskip

[M] Yu.~Manin. {\it Generating functions in algebraic geometry
and sums over trees.} In: The Moduli Space of Curves, ed. by
R\. Dijkgraaf, C\. Faber, G\. van der Geer, Progress in Math\.
vol\. 129, Birkh\"auser, 1995, 401--418.

\enddocument